\documentclass[11pt,twoside]{article}
\makeatletter
\newcommand{\specificthanks}[1]{\@fnsymbol{#1}}
\makeatother

\usepackage{latexsym}
\usepackage{amssymb, amsbsy, amsmath, amsfonts, amssymb, amscd, color}
\usepackage{graphicx}
\usepackage{float}
\usepackage[colorlinks=true]{hyperref}
\newcommand{\commentout}[1]{}
\newcommand{\R}{\mathbb{R}}

\newcommand {\al} {\alpha}

\newcommand {\Chi} {{\bf \raise 2pt \hbox{$\chi$}} }
\newcommand {\f}   {\frac}
\newcommand {\p}   {\partial}
\newcommand{\dis}{\displaystyle}
\newcommand {\proof} {\noindent {\bf Proof}. }
\newcommand{\beq}{\begin{equation}}
\newcommand{\eeq}{\end{equation}}
\newcommand{\bea} {\begin{array}{rl}}
\newcommand{\eea} {\end{array}}
\newcommand{\beac} {\begin{array}{c}}
\newcommand{\eeac} {\end{array}}
\newcommand{\bepa}{\left\{ \begin{array}{l}}
\newcommand{\eepa} {\end{array}\right.}
\newtheorem{theorem}{Theorem}[section]

\newtheorem{remark}[theorem]{Remark}
\newtheorem{proposition}[theorem]{Proposition}

\newcommand{\qed}{{ \hfill
		{\unskip\kern 6pt\penalty 500 \raise -2pt\hbox{\vrule\vbox to 6pt{\hrule width 6pt
					\vfill\hrule}\vrule} \par}   }}

\title{A model for cost efficient Workforce Organizational Dynamics and its optimization}
	\author{{Beno\^ \i t Perthame}\thanks{Sorbonne Universit\'es, UPMC Univ Paris 06, Laboratoire Jacques-Louis Lions  UMR CNRS 7598, UPD, Inria de Paris team Mamba, F75005 Paris, France} 
	\and{Edouard Ribes}\thanks{Email: edouard.augustin.ribes@gmail.com}
	\and{ {Karim Touahri}\thanks{Email: touahrikarim91@gmail.com} }
	 \and{ {Delphine Salort}\thanks{Sorbonne Universit\'es, UPMC, Laboratoire de Biologie Computationnelle et Quantitative UMR CNRS 7238, F75005 Paris, France}}
	 }
	
\begin{document}
\maketitle
\begin{abstract}
	This paper presents a workforce planning model scalable to an entire hierarchical organization. Its main objective is to design a cost optimal target which leverages flexible workforce solutions while ensuring an efficient promotional flux. The value of this paper lies in its proposal of an adequate flexibility rate using various solution types and in its discussion about external hiring ratios. The mathematical structures of the models are analyzed and numerical simulations illustrate the theoretical background.  
\end{abstract}
\noindent{\makebox[1in]\hrulefill}\newline
\section{Introduction}

To our knowledge, human resource planning has mainly focused on providing insights on hiring and training decisions.  This has boiled down to questions of workforce planning \& knowledge acquisition \cite{Bordoloi2001169, Burger_2016}. According to a recent review \cite{VandenBergh2013367}, these topics are becoming more and more workforce centric, as firms try to integrate employee preferences in their decisions. This naturally raised the question of career and employee lifecycle within the company and leads to the question of employee promotion strategy \cite{Klehe2011217, balancing,henderson}.
\\
The question of the employee lifecycle is complex, because it depends not only in the company purpose and market landscape but also in the individual employee dynamics. We have therefore investigated those dependencies within 3  complementary papers,  and this paper is the first one of the series.  More precisely, the main objective of this article is to present a  new partial differential equation model, perform its theoretical and numerical studies, in order to address the description and cost optimization of workforce lifecycle from a company standpoint.  In this perspective, we use the notion of workforce flexibility and notably the use of temporary contracts to yield cost efficiency, which is  standard among the operation research literature (see \cite{Qin201519} for a recent review).  However, in the present article, we  integrate a new  continuous formalism which allows to integrate mixed approaches in complex environment. The second paper focuses on turnover prediction and on identification of endogenous and exogenous parameters that influence voluntary turnover decisions (\cite{ribes:hal-01556746}). The turnover  rates are also used in the model data in section~\ref{sec:simu}. Finally the third paper, under preparation, adopts a mean field game approach and tries to reconnect company policies with individual aspirations. 
\\
\\
The article is organized as follows. In section \ref{sec:OD}, we build a dynamic organization representation based upon recent researches in this field \cite{balancing} by taking into account the time spent by personnel in a job qualification class.  While a simple model, which already use the idea of a waiting time before promotion has been used, with a simple fixed time delay model in \cite{henderson}, in our model, we integrate a more accurately control of this notion of waiting time, by adopting a continuous approach.
In  section \ref{sec:flex},  we determine  theoretically the best strategy to adopt in order to minimize the total labor cost. To this,  we assume that  the cost can be modulated via the recruitment of   temporary contracts. Indeed, the possible cost burden of having an organization which promotion waiting times are important will be the center of this investigation. More precisely, a cost optimal solution leveraging a flexible workforce will be proposed and fully characterized.
In section \ref{sec:simu}, we illustrate our theoretical results via numerical simulations. We first implement the organizational model with a strategy which allows to reduce the external hiring with specific rules. In a second part, we test, via genetic algorithms,  the  optimal strategy to minimize the labor cost.   We  show that a mix of "floaters" and "temporary workers"(see \cite{Wild199395} and \cite{Harper2010} for examples of the individuals approaches), can be used to efficiently design a full organizational footprint. Finally, in section \ref{sec:conc},  we finish our article by a conclusion an perspectives.

\paragraph{Important legal remarks} The findings and opinions expressed in this paper are those of the authors and do not reflect any positions from any company or institution. Finally, please bear in mind that to protect confidentiality numbers have been disguised in a way that preserves the same analysis and conclusions as the actual case study.

\section{Organization dynamics}	

We begin with  the main  assumptions and model we use in order to describe the dynamic of the workforce when employees may  be promoted on several hierarchical jobs. We discuss several situations with respect to the possibility of external hiring and we propose a theoretical study of this model, in order to determine the well-posedness of the related system of equations, the possible  stationary states and the necessity of external hiring, depending upon  parameter combinations, in order to  to obtain a well-posed problem.  This model is the basis of the study of the optimization cost.
\label{sec:OD}
\subsection{Organization model}
We represent a company by a set of $L$ jobs ordered in terms of level of responsibility from $j=1$ to $L$. We denote by $\rho_j(t,s)$ the number of workers which at time $t$, have seniority $s$ in their job class. 
We assume that:
\begin{itemize}
	\item in each job $j \in \{1,..,L\}$, the total number of workers in the job category $j$ across seniority categories is constant over time and equal to~$N_j>0$.
	\item in each job $j$, workers can leave the company. We assume that worker turnover follows an exponential law with parameter~$\mu_j$.
	\item workers can also stay in the same job category $j$ and wait to be promoted from job $j$ to $j+1$ .The promotion can only occur if they have spent a time superior to ~$\tau_j \geq 0$  in their job $j$ . The associated promotion rate is defined as ~$P_j$.
	\item the company can also hire externally in each job $j$ with a rate $h_j$. We denote by $h_j(t) N_j$ the number of external workers recruited externally for the job~$j$.
\end{itemize}
These assumptions can be formalized using the following balance law for workers in each job~$j=1,...,L$, which takes into account the seniority increase with time, as well as turnover and  promotion:
\begin{equation}
\partial_t \rho_j+ \partial_s \rho_j+ \mu_j \rho_j + \mathbb{I}_{s \geq \tau_j} P_j (t)\rho_j(t,s)=0, \qquad s\geq 0, \; t \geq 0,
\label{eq:1}
\end{equation}
The next boundary condition describes the flux of workers in class $j=1,...L$ that newly arrive either by external hiring or internal promotion:
\begin{equation} 
\rho_j(t,s=0)= h_j N_j + P_{j-1} (t)A_{j-1} (t), \quad \qquad (P_0=0, \quad P_L=0), 
\label{eq:2}
\end{equation}
with the notation
\begin{equation} 
A_j(t)= \int_{s= \tau_j}^{+\infty} \rho_j(t,s)ds.
\label{eq:3}
\end{equation}
Notice that for the first job class of job, only external hiring is possible. Thus the above condition holds with $P_0 = 0$. Also, for the last class, promotion is no longer possible and thus $P_L =0$. To simplify the setting, we assume that retirement is included in the average departure rate $\mu_j$.
The initial state is denoted by $\rho_j^0(s)> 0$ and solutions remain positive. We recall that this type of system, called renewal equations, is standard and well understood~\cite{Perthame_book, gurtin1974non, farkas2002stability}, and has previously been used for workforce  planning purposes ~\cite{Doumic2017}. 

As the number of workers in the job category $j$ is assumed constant, the following constraint appears:
\begin{equation}
N_j= \int_{0}^{+\infty} \rho_j(t,s)ds >0 \quad \hbox{is given}.
\label{eq:4}
\end{equation}
Integrating  equation~\eqref{eq:1}, we find a consistency condition on the boundary
\begin{equation}	
\rho_j(t,s=0)=  \mu_j N_j+ P_j A_j(t).
\label{eq:5}
\end{equation}
Combining the boundary condition \eqref{eq:2} and the constraint \eqref{eq:5} leads to $\forall  j=1,...,L$:
\begin{equation}\label{bc1}
\mu_j N_j+ P_j A_j(t)= h_j(t) N_j+ P_{j-1}A_{j-1} (t), \qquad \quad (P_0=P_L=0)
\end{equation}
At this stage, it has to be noticed that, because we have imposed the $L$ constraints that the $N_j$'s are constant, we also need $L$  Lagrange multipliers to fulfill the equality~\eqref{bc1}. These are to be chosen among the $2L$ parameters $P_j$ and $h_j$ and thus some degrees of freedom remain. 
In a first pass, a perfect internal labor market will be considered. This means that no external hiring will exist except at the first job layer $j=1$ and that the $L$ free parameters are $h_1$, $P_j$ for $j=1,...L-1$. In a second pass, a proportion of external hires among the various job layers will be introduced. It will be shown that it results from a simple generalization of the first case.
\subsection{Organization behavior - no external hiring}
Under the preliminary assumption of a perfect internal labor market, the company's recruitment policy prohibits external recruitment except for the first job layer. This means that $h_j=0$  for all $ j \geq 2$. Therefore, the constraint \eqref{bc1} becomes:
\begin{equation}\label{bc1case2}
\bepa
\mu_j N_j+ P_jA_j(t)=  P_{j-1}A_{j-1} (t), \qquad \hbox{for }  j = 2,...,L , 
\\[2mm]
\mu_1 N_1 +P_1 A_1(t)= h_1(t) N_1.
\eepa
\end{equation}
\paragraph{Promotional fluxes} 
Assuming that these relation can be solved backward (that is $A_j(t) >0$, for all $t \geq 0$ and $j \in \{1, ..., L-1\}$), the equation \eqref{bc1case2} imposes that for  the job level $j=L$:
$$
\mu_L N_L=  P_{L-1}A_{L-1}(t)
$$
This means that the promotion flux $ P_{L-1}A_{L-1}$ is independent of time. This can be further cascaded down the organization layers with $j=L-1$ to get:
$$
\mu_{L-1} N_{L-1}+ P_{L-1}A_{L-1}= P_{L-2}A_{L-2}, 
$$
and thus:
$$
\mu_{L-1} N_{L-1}+ \mu_L N_L =  P_{L-2}A_{L-2}.
$$ 
Iterating the same process shows that:
$$
\sum_{l=j}^{L} \mu_l N_l=  P_{j-1} A_{j-1}, \qquad  \forall j, \quad   2\leq j \leq L.
$$
Therefore, the promotion fraction $P_j(t)$  only depends upon $A_j(t)$ thanks to the relation:
\begin{equation}\label{without:Cj}
P_{j}(t)=  \frac{C^{no}_{j+1}}{A_{j}(t)}, \qquad C^{no}_j :=  \sum_{l=j}^{L} \mu_l N_l.
\end{equation}
Departing from equations  \eqref{eq:1} and   \eqref{eq:2}, this generates the following system of $L$ nonlinear equations:
\begin{equation}\label{eqf0}
\partial_t \rho_L + \partial_s \rho_L + \mu_L \rho_L=0, \qquad \rho_L(s=0,t)= \mu_L N_L,
\end{equation}
and for $1 \leq j \leq L-1$,
\begin{equation}\label{eqf}
\partial_t \rho_j + \partial_s \rho_j + \mu_j \rho +  \frac{ C^{no}_{j+1} }{A_{j}(t)} \mathbb{I}_{s \geq \tau_j}  \rho_j=0, \qquad \rho_j(s=0,t)=  \sum_{l=j}^{L} \mu_l N_l .
\end{equation}
The boundary condition for $j=1$ then states that $h_1N_1= \sum_{l=1}^{L} \mu_l N_l$. The interpretation is clear: the total number of hired workers at job level $j=1$ should compensate the total number of workers leaving the company. 
\paragraph{ Well-posedness of equation \eqref{eq:1}, assuming \eqref{bc1case2}}
Conditions on the initial data are now computed. This introduces a well posed problem for equation \eqref{eq:1}, when the parameters are chosen as in \eqref{bc1case2}.  This leads us to describe the set of initial data such that $A_j(t) >0$, for all $t \geq 0$ and $j \in \{1, ..., L-1\}$. The following proposition therefore holds:
\begin{proposition}\label{prAjnoh}
	Let $\rho^0_j(s) \in L^1(\R^+)$ such that for all $t \leq \tau_j$
	\begin{equation}\label{without:condin}
	e^{-t \mu_j}\int_{t}^{\tau_j} \rho_j^{0}(s-t) ds  +  \frac{1}{\mu_j} (1-e^{-\mu_j t})   \sum_{l=j}^{L} \mu_l N_l < N_j,
	\end{equation}
	Assuming that for all $j= 1,..., L-1$, the model parameters are such that:
	\begin{equation}\label{without:condParam}
	\overline{A_j} :=  \frac{\mu_j N_j e^{-\mu_j \tau_j}  -  (1-e^{- \mu_j \tau_j})  C^{no}_{j +1}}{\mu_j} >0.
	\end{equation}
	Then, equation \eqref{eqf} is well posed. We have $A_j(t)>0$ for all $t \geq 0$ and $A_j(t)=  \overline{A_j}$   for all $ t \geq \tau_j$.
\end{proposition}
\begin{remark}
	Let us mention that proposition \ref{prAjnoh} implies that equation \eqref{eqf} is linear as soon $t \geq \tau_j$ because $A_j$ then becomes constant.
\end{remark}
{\bf Proof of proposition \ref{prAjnoh}}
Estimating \eqref{without:condin} with $t=0$ is equivalent to choosing an initial data such that $A_j(0)>0$. Moreover, as long as $A_j(t)>0$, using the explicit formula for the solution of equation \eqref{eqf} via the method of characteristics, we find that for $t \leq \tau_j$,
$$\rho_j(t,s)= \rho_j^0(s-t) e^{-\mu_j t} \quad  \hbox{ if }     \tau_j \geq s\geq t$$
and 
$$\rho_j(t,s)=  \sum_{l=j+1}^{L} \mu_l N_l e^{-\mu_j s} \quad  \hbox{ if } s \leq t \leq \tau_j.$$
Using the mass conservation, we obtain that for $t \leq \tau_j$:
$$A_j(t)=  N_j - \int_{0}^{\tau_j} \rho_j(t,s)ds= N_j- e^{-t \mu_j}\int_{t}^{\tau_j} \rho_j^{0}(s-t) ds  -  \frac{1}{\mu_j} (1-e^{-\mu_j t})   \sum_{l=j}^{L} \mu_l N_l,$$
and so $A_j(t)>0$ because of condition \eqref{without:condin}.
Now, for $t \geq \tau_j$, with the method of characteristics, we obtain that :
$$\rho_j(t,s)=  \sum_{l=j}^{L} \mu_l N_l e^{-\mu_j s} \hbox{ if } s  \leq \tau_j$$
which implies that :
$$A_j=   \frac{\mu_j N_j e^{-\mu_j \tau_j}  -  (1-e^{- \mu_j \tau_j}) \sum_{l=j+1}^{L} \mu_l N_l}{\mu_j} , \quad \forall t \geq \tau_j.$$
\hfill $\square$

\paragraph{Stationary states and convergence of the solution of Equations   \eqref{eqf0}--\eqref{eqf}}
The following proposition expresses the fast convergence to a steady state
\begin{proposition}\label{propstatio}
	Assume that for all $j \in \{1,...,L-1 \}$, we have $\overline{A_j}>0$. Then, there exists a unique stationary state of Equations \eqref{eqf0}-\eqref{eqf} given by
	$$  
	\overline{\rho_L}(s)=  \mu_L N_L e^{- \mu_L s},
	$$
	$$  
	\overline{\rho_j}(s)=  e^{-\mu_j s + B_j(\tau_j-s)_-} \left( \sum_{l=j}^{L} \mu_l N_l   \right).
	$$
	Moreover, for initial data chosen as in proposition  \ref{prAjnoh}, there is a constant $C$ such that for all $j \in \{1, ..., L\}$ and $t \geq \max_{j} \tau_j$
	$$
	\int_{0}^{+\infty} | \rho_j(t,s)- \overline{ \rho_j}(s)|ds \leq C e^{-\mu_j t} \int_{0}^{+\infty} | \rho_j(0,s)- \overline{ \rho_j}(s)|ds.
	$$
\end{proposition}
{\bf Proof of proposition \ref{propstatio}}
Using that for $t \geq \tau_j$, $A_j(t)= \overline{A_j}$, the stationary state of equation \eqref{eq:1} is the stationary state of the following linear equation:
\begin{equation}\label{eqflin}
\partial_t \rho_j + \partial_s \rho_j + \mu_j \rho + B_j  \mathbb{I}_{s \geq \tau_j}  \rho_j=0, \qquad \rho_j(s=0,t)=  \sum_{l=j}^{L} \mu_l N_l
\end{equation} 
with $\displaystyle B_j:=\frac{ \sum_{l=j+1}^{L} \mu_l N_l}{\overline{A_{j}}}$ which leads to the formula of the stationary states given in Proposition \ref{propstatio}. To prove asymptotic convergence, we set $m_j=  \rho_j(t,s) - \overline{\rho_j}(s)$. We have  for all $j \in \{1, ...,L\}$ 
$$\partial_t |m_j| + \partial_s |m_j| + \mu_j |m_j| =0.$$
Integrating this equation and using the Gronwall lemma, we find the result of convergence which ends the proof of Proposition~\ref{propstatio}. \hfill $\square$

\subsection{Organization behavior - with external hiring}

We remove the constraint $h_i=0$ and allow the external organization to source itself externally. We adapt the method developed in the previous subsection to this situation with external sourcing and we show that, with a more general condition on parameters, it yields similar results.
Assuming again that for all $j$, $A_j(t)>0$, we introduce the hiring ratio  $\alpha_j >1$, for job class $j$ according to the below equation:
\beq
\mu_j N_j+ P_jA_j(t)=P_{j-1}A_{j-1} (t)+h_j.N_j=  \alpha_j \;   P_{j-1}A_{j-1} (t), \quad 2\leq j \leq L,
\label{ei:1}
\eeq
This gives us the external hiring rate $h_j(t)$ through the relation:
\begin{equation}\label{eqalpha}
\alpha_j= 1+\frac{h_j(t)N_j(t)}{P_{j-1}A_{j-1} (t)}.
\end{equation}
In other words, $\alpha_j $ is an interesting parameter for the establishment of human resources policies as $\alpha_j -1 $ represents the ratio of external over internal hires at level $j$. 
\paragraph{Promotional fluxes}
Following the argument in the previous subsection, we infer from \eqref{ei:1}, successively for $j=L$ and for $j=L-1$, 
$$
\mu_L N_L= \alpha_L  P_{L-1}A_{L-1} ,
$$ 
$$
\mu_{L-1}N_{L-1} + \alpha_L^{-1} \mu_L N_L= \alpha_{L-1} P_{L-2}A_{L-2} .
$$
Iterating this process for $2 \leq j \leq L-1$, leads to
$$
\mu_j N_j + \sum_{l=j+1}^{L}\mu_l N_l .(\prod_{k=j+1}^{l}\frac{1}{\alpha_k})=   \alpha_{j} P_{j-1} A_{j-1}, 
$$
that is also written
$$
C_{j}:= \sum_{l=j}^{L}\mu_l N_l .(\prod_{k=j}^{l}\frac{1}{\alpha_k})=   P_{j-1} A_{j-1}, 
$$
or also 
$$
P_j(t)=   \frac{ C_{j+1}}{ A_j(t)} , \qquad j=1,..., L-1.
$$
This generates the following systems of $L$ nonlinear equations:
\begin{equation}\label{eqf0bis}
\partial_t \rho_L + \partial_s \rho_L + \mu_L \rho_L=0, \qquad \rho_L(s=0,t)= \mu_L N_L ,
\end{equation}
and for $1 \leq j \leq L-1$,
\begin{equation}\label{eqfbis}
\partial_t \rho_j + \partial_s \rho_j + \mu_j \rho +  \frac{C_{j+1}}{A_{j}(t)} \mathbb{I}_{s \geq \tau_j}  \rho_j=0, \qquad \rho_j(s=0,t)= \mu_jN_j+C_{j+1} \;(= \al_j C_j).
\end{equation}
\paragraph{Study of equations  \eqref{eqf0bis}--\eqref{eqfbis}} 
The structure of equations \eqref{eqf0bis}--\eqref{eqfbis} are exactly the same as equations \eqref{eqf0}--\eqref{eqf}. Hence, in this section, we just state our results without proofs since they follow exactly the preceding section. 
\newline 
Given $\alpha_j \geq 1$, we begin with conditions on the initial data to have a well posed problem for equations \eqref{eqf0bis}--\eqref{eqfbis}. The following proposition holds:
\begin{proposition}\label{propinbis}
	Let $\rho^0_j \in L^1(\R^+)$, $1\leq j \leq L$, be such that for all $t \leq \tau_j$
	\begin{equation}\label{with:condin}
	e^{-t \mu_j}\int_{t}^{\tau_j} \rho_j^{0}(s-t) ds  +  \frac{1}{\mu_j} (1-e^{-\mu_j t}) ( \mu_jN_j+C_{j+1} ) < N_j ,
	\end{equation}
	Assume  that for all $j= 1,..., L-1$, the model parameters are such that:
	\begin{equation}\label{with:condParam}
	\widetilde{A_j} := \frac{\mu_j N_j e^{-\mu_j \tau_j}  -  (1-e^{- \mu_j \tau_j}) C_{j+1}}{\mu_j} >0 .
	\end{equation}
	Then equation \eqref{eqf} is well posed. We have $A_j(t)>0$ for all $t \geq 0$ and 
	$$
	A_j(t)=  \widetilde{A_j}  , \quad \forall t \geq \tau_j.
	$$
\end{proposition}

\begin{remark} Let us mention that, given $(\tau_j)_{j \in \{1,..,L-1\}}$, $(\mu_j)_{j \in \{1,..,L\}}$ and $(N_j)_{j \in \{1,..,L\}}$, we can always find $(\alpha_j)_{j \in \{1,..,L\}}$ such that for all $ j \in \{1,..,L\}$, $\widetilde{A_j} >0$. In the same way, given an initial data in $L^1$ such that for all  $j \in \{1,..,L\}$, $A_j(0)>0$, we can always find $(\alpha_j)_{j \in \{1,..,L\}}$  such that  $A_j(t)>0$ for all $t \geq 0$.
\end{remark}

\begin{proposition}\label{propstatiobis}
	Assume that for all $j \in \{1,...,L-1 \}$, we have $\widetilde{A_j}>0$. Then, there exists a unique stationary state of equations \eqref{eqf0bis}-\eqref{eqfbis} given by:
	$$ 
	\widetilde{\rho_L}(s)=  \mu_L N_L e^{- \mu_L s},
	$$
	$$ 
	\widetilde{\rho_j}(s)=  e^{-\mu_j s - \frac{C_{j+1}}{\widetilde{A_j}}(s- \tau_j)_+} \left( \mu_j N_j+C_{j+1}   \right).
	$$
	Moreover, for initial data chosen as in Proposition~\ref{propinbis}, there is a constant $C$ such that for all $j \in \{1, ..., L\}$ and $t \geq \max_{j} \tau_j$
	$$ 
	\int_{0}^{+\infty} | \rho_j(t,s)- \widetilde{ \rho_j}(s)|ds \leq C e^{-\mu_j t} \int_{0}^{+\infty} | \rho_j(0,s)- \widetilde{ \rho_j}(s)|ds.
	$$
\end{proposition}

\paragraph{Minimal external hiring} A company may prefer to promote internally to capitalize on its own workforce and provide its employees with better careers opportunities to minimize their turnover rate. To do so, the procedure is to choose the `minimal values' of the vector $(\al_1,...,\al_L)$ so as to impose the constraint~\eqref{with:condParam}. This can be performed by a descending algorithm. 
Departing from $\widetilde{A_{L-1}}>0$, we find :
$$
C_{L} := \f{\mu_L N_L}{\al_L}  <  \frac{\mu_{L-1} N_{L-1} e^{-\mu_{L-1} \tau_{L-1}}}{(1-e^{- \mu_{L-1} \tau_{L-1}}) } ,
$$
$$
\al_L^{min} = \max \left(1, \frac{(1-e^{- \mu_{L-1} \tau_{L-1}})\mu_L N_L }{\mu_{L-1} N_{L-1} e^{-\mu_{L-1} \tau_{L-1}}} \right).
$$
This allows to compute $\al_{L-1}^{min}$ because we impose  $\widetilde{A_{L-2}}>0$, which is 
$$
C_{L-1} := \f{1}{\al_{L-1}} \left[  {\mu_{L-1} N_{L-1}}  + \f{\mu_L N_L }{\al_L} \right] <  \frac{\mu_{L-2} N_{L-2} e^{-\mu_{L-2} \tau_{L-2}}}{(1-e^{- \mu_{L-2} \tau_{L-2}}) } ,
$$
Hence the formula for 
$$
\al_{L-1} ^{min} =\max\left(1,  \frac{ 1-e^{- \mu_{L-2} \tau_{L-2}} }{\mu_{L-2} N_{L-2} e^{-\mu_{L-2} \tau_{L-2}}}
\big[  {\mu_{L-1} N_{L-1}}  + \f{\mu_L N_L }{\al_L^{min}} \big] \right).
$$

\section{Organization cost structure and flexible workforce}
\label{sec:flex}

In the previous section, we have shown that it is possible to control the dynamic of an entire organization in a continuous framework. The aim of   this section  is   to deliver a theoretical setting to understand how to make the labor cost efficient and optimal. The problematic being that, for a given position, permanent workers in an organization have an increasing salary with seniority. As such, it may prove detrimental to have them waiting for a promotion. Flexible workforce solutions will therefore be discussed to minimize the overall labor cost burden of the organization. To this, we first study  the impact of the use of temporary workers in the labor cost and then integrate a mixed approach of "floaters" and temporary workers.

\subsection{Leverage temporary workers to speed up the organization}
\label{sec:temporary}

\paragraph{General problem statement} Using the notations in section~\ref{sec:OD}, we consider that  the population $N_j$ in each job $j \in {1,...,L}$ can be divided in two subcategories:
\begin{itemize}
	\item The permanent population $N_j^p >0$ which evolves within the organization with the same dynamics as described earlier and which compete for promotion. The permanent population distribution across seniority level at steady state will be referred to as $\rho_j^{p}(s)$.
	\item The temporary population $N_j^t >0$ which is used by the organization as a buffer to meet its workload requirements and will not be used for promotion.
\end{itemize}
We denote by  $p_j \in [0,1]$ the proportion of the population $N_j$ in level $j$ under a permanent contract. We may write
\begin{equation}
N_j=N_j^p+N_j^t=(1-p_j) N_j+ p_j N_j .
\end{equation}
The organization cost structure is supposed to obey the following rules at each level $j$:
\begin{itemize}
	\item the cost of each position $w_j^p(s)>0$ in the permanent population $N_j^p$ is growing with seniority $s$,
	\item the cost of each temporary contract $w_j^t$ is assumed to be constant. We also assume that $w_j^t>w_j^p(0)$ to account that a temporary workforce may come at a premium over a freshly hired permanent employee.
\end{itemize}
The overall operating cost of the organization $Cost_{org}$  is therefore defined as the sum of the operating costs $Cost_j$ at each level $j$, 
\begin{equation}
Cost_{org}=\sum_{j=1}^{L}Cost_j
\end{equation}
with
\begin{equation}
Cost_j= (1-p_j) N_j \; w_j^t+\int_0^\infty \rho_j^{p}(s) w_j^p(s) ds.
\end{equation}
\paragraph{Computing the organization cost}
According to Proposition~\ref{propstatiobis} and Equation~\eqref{eqfbis}, we have 
\begin{itemize}
	\item $\rho_L^p(s)=  \mu_L N_L\; p_L.e^{- \mu_L s}$,
	\item  $\rho_j^p(s)= (\mu_j N_j+C_{j+1})  e^{-\mu_j s + \frac{C_{j+1}}{{\widetilde A}_j}(\tau_j-s)_-}  $ for  $1\leq j<L$,
\end{itemize}
where
$$
C_{j}(p_{j},...,p_L): =\dis \sum_{l=j}^{L}\; \mu_l N_l.p_l.\prod_{k=j}^{l}\frac{1}{\alpha_k},
$$
$$
{\widetilde A}_j= \dis \frac{\mu_j N_j.p_j. e^{-\mu_j \tau_j}  -  (1-e^{- \mu_j \tau_j}) C_{j+1}}{\mu_j}.
$$
Therefore, we may define the optimal cost  organization by:
\begin{equation}\label{costoptimA}
Cost_{org}^{opt}= \min_{ (p_1,...,p_L)\in [0,1]^L} \sum_{j=1}^{L} \left\{(1-p_j) N_j \; w_j^t+\int_0^\infty \rho_j^{p}(s) w_j^p(s) ds  \right\}.
\end{equation}
Note that, similarly to the previous condition~\eqref{with:condParam},  we need to impose a limitation of the $p_j$ in order to satisfy ${\widetilde A}_j >0$, namely 
\beq
p_j >  \f{ (1-e^{- \mu_j \tau_j}) C_{j+1}}{\mu_j N_j. e^{-\mu_j \tau_j} } =: p_j^{\rm min}
\label{temporary:condin}
\eeq
where $p_j^{\rm min}$ depends only upon $p_{j+1},...,p_L$.
\paragraph{Solution Approach} In order to derive explicit formulae, we now assume that there is no external hiring ($\alpha_j \equiv 1$) and that the cost of permanent employees grows exponentially with a limited rate $r$ (compared to attrition)  
$$
w_j^p(s)  = w_j^0 \; e^{rs}, \qquad r < \mu_j , 
$$
Then, we can obtain explicitly the cost at level $j$ with formula which can be implemented within a standard optimization algorithm, see Section~\ref{sec:numFlex}.
\begin{proposition} The following expressions holds for the costs $Cost_j$, 
	\beq
	Cost_j= (1-p_j) N_j \; w_j^t +\f{w_j^0  (\mu_j N_j p_j+C_{j+1})  }{\mu_j-r} \left[  1 -  \f{  \mu_j C_{j+1} e^{r \tau_j} }{(\mu_j -r )\alpha_j C_j +r C_{j+1} e^{\mu_j \tau_j} } \right], \qquad 1 \leq j<L ,
	\label{flex:formulaj}
	\eeq 
	\beq
	Cost_L= (1-p_L) N_L \; w_L^t +\f{w_L^0 \alpha_L C_L}{\mu_L-r} .
	\label{flex:formulaL}
	\eeq 
\end{proposition}
\proof These formulas are derived as follows for $j<L$. Note that we do not copy the calculation for the case $j=L$. Departing from the formula for $\rho_j^p(s)$  in Section~\ref{sec:flex}, 
$$\bea
\f1{w^0_j} \dis\int_0^\infty \rho_j^{p}(s) w_j^p(s) ds& =  (\mu_j N_j p_j+C_{j+1})  \dis \int_0^\infty e^{-\mu_j s+rs +(r+ \frac{C_{j+1}}{A_j}) (\tau_j-s)_-} ds
\\[10pt]
& =  (\mu_j N_j p_j+C_{j+1})  [ \f{1- e^{-(\mu_j-r) \tau_j}}{\mu_j-r} +  \frac{ e^{-(\mu_j -r)\tau_j} } {\mu_j - r + \frac{C_{j+1}}{A_j} }]
\eea$$
Therefore, since $ \mu_j {\widetilde A}_j= e^{-\mu_j \tau_j}\alpha_j C_j- C_{j+1}$, we have 
$$\bea
\dis  \dis\f{\mu_j-r}{w^0_j}\int_0^\infty \rho_j^{p}(s) w_j^p(s) ds& =  (\mu_j N_j p_j+C_{j+1})  \left[ 1- e^{-(\mu_j-r) \tau_j}+  \frac{ (\mu_j-r) {\widetilde A}_j e^{-(\mu_j-r)\tau_j} } {e^{-\mu_j \tau_j} \alpha_jC_j - r{{\widetilde A}_j} } \right]
\\[10pt]
& = (\mu_j N_j p_j +C_{j+1})   \left[  1 + e^{-(\mu_j-r) \tau_j} \; \f{ -  ( e^{-\mu_j \tau_j} \alpha_jC_j - r {\widetilde A}_j ) +  (\mu_j -r) {\widetilde A}_j }{ e^{-\mu_j \tau_j} \alpha_jC_j - r {\widetilde A}_j  } \right]
\\[10pt]
& = (\mu_j N_j p_j+C_{j+1})  \left[  1 - e^{-(\mu_j-r) \tau_j} \; \f{  C_{j+1}  }{ e^{-\mu_j \tau_j}\alpha_j C_j - r {\widetilde A}_j  } \right].
\eea$$
We end up with the expression:
$$
\f{\mu_j-r}{w^0_j}\int_0^\infty \rho_j^{p}(s) w_j^p(s) ds= 
(\mu_j N_j p_j+C_{j+1})   \left[  1 -  \f{  \mu_j C_{j+1} e^{r \tau_j} }{\mu_j\alpha_j C_j-r\alpha_jC_j +r C_{j+1} e^{\mu_j \tau_j} }
\right]. 
$$
This formula immediately yields \eqref {flex:formulaj}. 
\qed

\paragraph{Optimality conditions}  The specific form of the cost in formula~\eqref{flex:formulaL} has an {\em ascending} property because $Cost_j$ only depends on $p_j,...p_L$. This allows us to write: 
$$
Cost_{org}^{opt}= \min_{p_L} \Big[ Cost_L + \min_{p_{L-1}} \big[ Cost_{L-1}+  \min_{p_{L-2}}[ Cost_{L-2} ... +\min_{p_1} Cost_1  ...]\big] \Big].
$$
Notice that this structure is conflicting with the {\em descending} structure for the minimal values of the $\al_j$. For this reason analytical formula for a global minimizer in $(p_j, \al_j)_{j=1,...L}$ is not possible. 
We exemplify the choices of optimal proportions of permanent workers $p_j$'s in the next paragraph, assuming the $\al_j$'s are known.
\paragraph{A specific case of optimality conditions} 
Computing the optimal conditions remains rapidly very tedious, at least theoretically. But rather than providing directly numerical simulation results (that are in any cases developed in the end of this paper), we will restrict our strategy of leverage of temporary workers only in the first job or the two first jobs to derive and discuss closed formulas. 

\vspace{0,3cm}
\noindent \underline{Case 1: Allow contractors only in the lowest hierachical level.}
In this setting, the cost function depends only on the variable $p_1$ because $p_2=p_3=...=p_L=1.$.   Because only $C_1$, among the $C_j$'s, depends on $p_1$, from \eqref{flex:formulaj},  we then have to tackle the minimums with respect to the variable $p_1$ of the function 
$$
 Cost_1 (r,p_1)=  (1-p_1) N_1 w_1^t+ \f{w_1^0  \alpha_1 C_1 }{\mu_1-r} \left[  1 -  \f{  \mu_1 C_{2} e^{r \tau_1} }{(\mu_1 -r )  \alpha_1 C_1 +r C_{2} e^{\mu_1 \tau_1} } \right].
 $$
We have 
$$
 Cost_1 '(r,p_1)=  N_1 \left( -w_1^t+ \frac{w_1^0 \mu_1 }{\mu_1-r}  \big(1 -  \frac{r \mu_1(C_2)^2 e^{(r+ \mu_1) \tau_1} }{((\mu_1 -r ) \alpha_1 C_1 +r C_{2} e^{\mu_1 \tau_1}  )^2  }\big) \right)\cdotp
$$ 
We deduce that 
$$
 Cost_1 ''(r,p_1) =  \frac{2  N_1 w_1^0 r \mu_1^2 (C_2)^2 e^{(r+ \mu_1) \tau_1}}{((\mu_1 -r ) \alpha_1 C_1 +r C_{2} e^{\mu_1 \tau_1}  )^3} >0  \quad \hbox{ on } ]0,1[.
$$
Hence, $ Cost_1 $ is a strictly convex function with respect to the variable $p_1$ and there exists a unique minimum $p_1^*$ of $ Cost_1 $ on $[0,1]$. Moreover, if  $p_1^* \in ]0,1[$, we have $ Cost_1'(r,p_1^*)=0$. 
We can further analyze the optimal choice of $p_1$.The optimality condition $ Cost_1 '(r,p_1)=0$ gives us :
$$
1- \f{\mu_1-r} {\mu_1} \; \f{w_1^t }{w^0_1} =  \frac{r \mu_1(C_2)^2 e^{(r+ \mu_1) \tau_1} }{\big((\mu_1 -r ) p_1\mu_1 N_1  + C_{2}(r e^{\mu_1 \tau_1} + \mu_1-r ) \big)^2  }
 < \frac{r }{\mu_1} e^{(r- \mu_1) \tau_1}
$$
because we have assumed from \eqref{temporary:condin} that  $p_1 \mu_1 N_1 > (e^{\mu_1 \tau_1}-1)C_2$ (a relation written as $p_1>p_1^{\rm min}$). Therefore:
\\
$\bullet$ when $w_1^t$ is low the optimal cost is achieved with the smaller tenure rate $p_1^*=p_1^{\rm min}$ which allows to ensure internal promotion for level 2 (see \eqref{temporary:condin}),
\\
$\bullet$ when $w_1^t$ is large, the minimal cost is also reached by saturating the constraint, but now on the other extreme $p_1^*=1$ which means that there is no temporary workers,
\\
$\bullet$ the middle range where the optimal hiring rate satisfies $p_1^{\rm min}< p_1^*<1$ occurs when the parameters are such that:
$$
\frac{r \mu_1(C_2)^2 e^{(r+ \mu_1) \tau_1} }{\big((\mu_1 -r ) \mu_1 N_1  + C_{2}(r e^{\mu_1 \tau_1} + \mu_1-r ) \big)^2  } < 1- \f{\mu_1-r} {\mu_1} \; \f{w_1^t }{w^0_1}   < \frac{r }{\mu_1} e^{(r- \mu_1) \tau_1}.
$$
\vspace{0,3cm}
\noindent \underline{Case 2: Allow contractors only in the two first levels.}  In this setting, we assume that  $p_3=...=p_L=1.$ This case is much more complex to deal with, for instance,  conversely to the previous case,  convexity of the cost function $cost_{org}$ is lost, even in the simpler situation where $L=2$.
Hence, here, we restrict our analysis to the derivation of the relations leading to a critical point, that is $\nabla Cost_{org}=0.$ 
\\[2mm]
As before, because only $C_1$, among the $C_j$'s, depends on $p_1$, from \eqref{flex:formulaj}, we obtain the first order optimality  condition $\f{\p Cost_{org}}{\p p_1} =0 $ with the relation:
\beq
(\mu_1 -r ) N_1 w_1^t =w_1^0 \mu_1 N_1  \left[ 1 -  \frac{r \mu_1(C_2)^2 e^{(r+ \mu_1) \tau_1} }{((\mu_1 -r ) \alpha_1 C_1 +r C_{2} e^{\mu_1 \tau_1}  )^2  } \right] .
\label{temp:2_1}
\eeq
The next first order optimality  condition $\f{\p Cost_{org}}{\p p_2} =0 $ is written
$$\bea
 N_2 w_2^t = & \mu_2 N_2  \left\{ \f{w_1^0}{\mu_1-r}  \left[  1 -  \f{  \mu_1 (C_{2}+\alpha_1 C_1) e^{r \tau_1} }{(\mu_1 -r ) \alpha_1  C_1 +r C_{2} e^{\mu_1 \tau_1} } 
  +  \alpha_1 C_1  \f{   \mu_1 C_{2} e^{r \tau_1} \big( (\mu_1 -r ) +r e^{\mu_1 \tau_1} \big)} { \big( (\mu_1 -r ) \alpha_1 C_1 +r C_{2} e^{\mu_1 \tau_1}\big)^2}   \right] \right.
 \\[15pt]
 & \left. +
 \f{w_2^0}{\mu_2-r}  \left[  1 -  \f{  \mu_2 C_{3} e^{r \tau_2} }{(\mu_2 -r ) \alpha_2 C_2 +r C_{3} e^{\mu_2 \tau_2} } +  \alpha_2 C_2  \f{(\mu_2 -r )  \mu_2 C_{3} e^{r \tau_2}}{ \big( (\mu_2 -r )  \alpha_2 C_2 +r C_{3} e^{\mu_2 \tau_2}\big)^2}\right] \right\},
\eea$$
which,using the first condition,  can be simplified  as
\beq\bea
 N_2 w_2^t = & \mu_2 N_2  \left\{\f{w_1^t }{\mu_1}+  \f{w_1^0}{\mu_1-r}  \left[  -  \f{  \mu_1  \alpha_1 C_1 e^{r \tau_1} }{(\mu_1 -r ) \alpha_1 C_1 +r C_{2} e^{\mu_1 \tau_1} } 
  + \alpha_1 C_1  \f{   \mu_1 C_{2} e^{r \tau_1} r e^{\mu_1 \tau_1}} { \big( (\mu_1 -r ) \alpha_1 C_1 +r C_{2} e^{\mu_1 \tau_1}\big)^2}   \right] \right.
 \\[15pt]
 & \left. +
 \f{w_2^0}{\mu_2-r}  \left[  1 -  \f{  \mu_2 C_{3} e^{r \tau_2} }{(\mu_2 -r ) \alpha_2 C_2 +r C_{3} e^{\mu_2 \tau_2} } +  \alpha_2 C_2  \f{(\mu_2 -r )   \mu_2 C_{3} e^{r \tau_2}}{ \big( (\mu_2 -r ) \alpha_2C_2 +r C_{3} e^{\mu_2 \tau_2}\big)^2}\right] \right\}.
\eea
\label{temp:2_2}
\eeq

These two relations \eqref{temp:2_1}, \eqref{temp:2_2} define the possible critical points $p_1$, $p_2$ of the cost function in this case.

\subsection{Add floaters to account for task specificities}

\begin{figure}[!h]
\centering
\includegraphics[height=5cm]{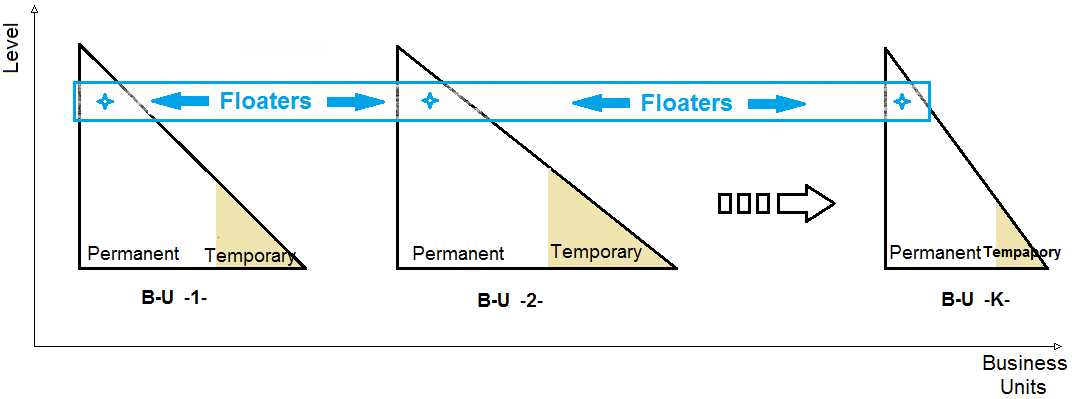}
\\[-4mm]
\caption{Representation of a company structure with temporay employees and floaters.}
\label{fig:BU}
\end{figure}

In the previous subsection, it was shown that the use of a temporary workforce can be optimized across all organizational levels to yield more cost efficiency by reducing the organization idleness and its associated costs. We now complete our model by adding a final layer to achieve what we expect to be a realistic description of an entire organization.
\paragraph{Final problem statement}Building on all the previous sections, assume that a company is actually the sum of $K$ business units. In each business unit $k \in [1,K]$, the workforce is organized around the same level structure $j \in [1,L]$ and we assume it is now structured in 3 parts (see Figure~\ref{fig:BU})
\begin{itemize}
	\item a $k$ business unit specific permanent population in level $j$ noted $N_j^{p,k}$ which follows the same dynamics as described earlier and which compete for promotion within the business unit.
	\item a temporary population $N_j^{t,k}$ which is used by the business units as buffer to meet its workload requirement and is not used for promotion.
	\item a permanent "floater" population $N_j^{float}$ that is not business unit specific and does not compete for promotion.
\end{itemize}
The introduction of "floaters" is of interest because it represents another type of workforce, where workers tend to be generalists that ripe the benefits of an organizational knowledge without competing for promotion. It also illustrates how an organization can benefit from transversal moves on top of vertical ones.
\\
Noting $p_j^k \in[0,1]$ the proportion of the population $N_j^k$ at level $j$ in business unit $k$ that is under a permanent contract and $g_j^k \in[0,1]$ the proportion of floaters, the following arises:
\begin{equation}
N_j^k=\underbrace{N_j^{k}p_j^k}_{\rm permanent} +\underbrace{ N_j^k.g_j^k}_{\rm floaters} +\underbrace{ N_j^k.(1-p_j^k -g_j^k)}_{\rm temporary}, \quad p_j^k+g_j^k \leq 1,  \qquad N_j^{p,k} := N_j^{k}p_j^k,
\end{equation}
and we may define for $ j \in [1,L]$,
\begin{equation}
N_j^{float}: =\sum_{k=1}^{K} N_j^k.g_j^k.
\end{equation}
From a cost standpoint, the following rules at each level $j$ and business unit $k$ apply:
\begin{itemize}
	\item the cost of each position $w_j^{p,k}(s)>0$ in the permanent population $N_j^{p,k}$ is growing with seniority~$s$,
	\item the cost of each temporary contract $w_j^{t,k}$ is assumed to be constant. We also assume that $w_j^{t,k}>w_{j}^{p,k}(0)$ to account that a temporary workforce may come at a premium over a freshly hired permanent employee,
	\item the cost of the "floaters" $w_j^{float}(s)>0$ is growing with seniority $s$.
\end{itemize}
This leads to a cost optimization program that is similar to~\eqref{costoptimA}, 
\begin{equation}\label{costoptimB}
\bea
Cost_{org}^{opt}= \! \! \min_{\footnotesize \beac p_j^k+g_j^k \leq 1, \\ k\in [1,K], j\in[1,L] \eeac}\! \dis \sum_{k=1}^{K}\sum_{j=1}^{L} \Big\{ & N_j^k.(1-p_j^k - g_j^k ).w_j^{t,k}+ \int N_j^k.g_j^k.w_j^{float}(s)e^{-\mu_j s} ds
\\
&+\dis\int w_j^{p,k}.\rho_j^{p,k}(s)ds \Big\}.
\eea
\end{equation}
\paragraph{Solution Aapproach}  We may analyze this cost functional and introduce the floaters averaged cost $w_j^{fa} = \int w_j^{float}(s)e^{-\mu_j s} ds$ and we obtain from \eqref{costoptimB}
$$
Cost_{org}^{opt}=\! \! \! \! \min_{\footnotesize \beac p_j^k+g_j^k \leq 1, \\ k\in [1,K], j\in[1,L] \eeac}\! \! \sum_{k=1}^{K}\sum_{j=1}^{L} 
\Big\{N_j^k.\big[(1-p_j^k  - g_j^k ).w_j^{t,k}+ .g_j^k.w_j^{fa}\big]+\int w_j^{p,k}.\rho_j^{p,k}(s)ds \Big\}.
$$
We first minimize over $g_j^k $ with the constraint $0 \leq g_j^k  \leq 1-p_j^k $. Because the functional depends linearly on $g_j^k $ it is attained either with $g_j^k =0$ or $g_j^k =1-p_j^k $ depending on the values of $w_j^{t,k}$ and $w_j^{fa}$.
We find the simplified formula:
$$
Cost_{org}^{opt}=\! \! \! \! \min_{\footnotesize \beac p_j^k \leq 1, \\ k\in [1,K], j\in[1,L] \eeac}\! \! \sum_{k=1}^{K}\sum_{j=1}^{L} 
N_j^k.\big[(1-p_j^k). \min(w_j^{t,k},.w_j^{fa}) + \int w_j^{p,k}.\rho_j^{p,k}(s)ds \Big\}.
$$
We are back to a problem, similar to the case of temporary workers in section \ref{sec:temporary},  with the additional complexity of the business units.

\section{Algorithms,  simulations \& results}
\label{sec:simu}

It is possible to determine numerically the solution of the organizational model~\eqref{eq:1}--\eqref{eq:2} with the size constraint~\eqref{bc1} for each job class. It allows either to minimize external hiring or to implement practically the above strategies. Either we determine the highest possible internal promotion rates $P_j(t)$, $j=1,...,L-1$ and, when there is not enough workers ready to be promoted, minimal external hiring is used to compensate the gap. Or, we impose a proportion of external hiring. 
The method we use is to solve the system~\eqref{eq:1}--\eqref{eq:2}, after re-writing it under the form
\begin{equation}\label{nc1}
\bepa
\partial_t \rho_j+ \partial_s \rho_j+ \mu_j \rho_j + P_j (t)\rho_j(t,s)=  \mathbb{I}_{s < \tau_j} P_j (t) \rho_j(t,s), 
\\[2mm]
\rho_j(t,s=0) =  (\mu_j +P_j) N_j+ P_j B_j(t), \qquad  B_j(t) = P_j (t) \int_0^{\tau_j} \rho_j(t,s)ds,
\eepa	 
\end{equation}
with the coefficients $P_j$ and $h_j$ as in \eqref{bc1}.
\begin{table}[p] \centering
	\begin{tabular}{cccccc}
		\hline\hline
		$j  $  & 1   &2     & 3    & 4  &5  \\
		\hline
		$N_j$ & 5500 & 5200 &3800 &1800 & 500 \\
		$\mu_j$& 0.8  &  0.8 &  0.8 &  0.8 & 5. \\
		$\tau_j$&4  &4 & 4 &4 &4  \\
		$P_j$& 4.42  &0.52 & 0.28 & 0.54 & --  \\
		$h_j$& 0.28  &0. & 0. & 0. & 0.  \\
		$T_j$& 0.33 & 1.78 & 2.9 & 1.72 & 2.1 \\
		$RP_j$& 0.045  & 0.26& 0.38 & 0.26 & 0.14
		\\
		\hline\hline
	\end{tabular}
	\\[-2mm]
	\caption{(Low turnover) The promotion rates $P_j$, ready to be promoted ratio $RP_j$ and waiting for promotion times $T_j$ for a company with 5 class jobs of size $N_j$ (in thousands of employees), attrition rate $\mu_j= 0.08$ per year, and waited to be promoted time $\tau_j$ in years. Internal hiring is sufficient.}
	\label{table1} \end{table}	
\subsection{Optimal organizational model (external hirings)}
The advantage of this form comes from the explicit-implicit discrete version which is both stable and easy to implement:
\begin{equation}\label{ncd1}
\bepa
\f{\rho^{k+1}_{j,i} - \rho^{k}_{j,i}}{\Delta t }+ \f{\rho^{k}_{j,i} - \rho^{k}_{j,i-1}}{\Delta s} + (\mu_j +P^k_j ) \; \rho^{k+1}_{j,i}  = P_j^k \; \mathbb{I}_{s_i < \tau_j} \;  \rho^{k}_{j,i}, \qquad i \geq 1,
\\[2mm]
\rho^{k}_{j,0} =  (\mu_j +P^k_j) N_j + \Delta s \; P^k_j  \; \sum_{0<s_i \leq \tau_j} \rho^{k}_{j,i}.
\eepa	 
\end{equation}
Here we have used a time step $\Delta t$, a s-step $\Delta s$ and the superscript refers to time $t^k= k \Delta t$, the subscript $i$ refers to the seniority $s_i=i \Delta s$. The only stability condition comes from the CFL restriction on $\Delta t$, namely $\Delta t \leq \Delta s$. We have used here a standard upwind scheme for discretizing the $s$ derivative~\cite{bouchut2004nonlinear, leveque1992numerical}.
Notice that, assuming by iteration the constraint,
\begin{equation}\label{ncconstraint}
\Delta s \sum_{i} \rho^{k}_{j,i} =N_j, 
\end{equation}
we find:
$$
\Delta s \sum_{i\geq 1}\rho^{k+1}_{j,i} - N_j + \Delta t  (\mu_j +P^k_j )  \Delta s \sum_{i \geq 1}\rho^{k+1}_{j,i}  = \Delta t  \rho^{k}_{j,0} +  \Delta t \Delta s \; P^k_j  \; \sum_{0 <s_i \leq \tau_j} \rho^{k}_{j,i}.
$$
This also gives:
$$
\Delta s \sum_{i\geq 1}\rho^{k+1}_{j,i} + \Delta t  (\mu_j +P^k_j )  \Delta s \sum_{i \geq 1}\rho^{k+1}_{j,i}  = N_j + \Delta t  (\mu_j +P^k_j) N_j ,
$$
the solution of which is $\Delta s \sum_{i\geq 1}\rho^{k+1}_{j,i}=N_j$. This means that the discretization \eqref{ncd1} preserves the discrete version \eqref {ncconstraint} of the constraint~\eqref{bc1}. 

It remains to determine the hiring rates $h^{k+1}_j, \; P^{k+1}_j$. To do so, we argue backward on the job index~$j$. After computing $\rho^{k+1}_{j,i}$, and departing from~$j=L$ and $P^{k+1}_L=0$, we impose until $j=1$ with $P_0^k=0$, 
\begin{equation}\label{ncchoice}
h_{j}^{k+1} N_j+ P_{j-1}^{k+1} A^{k+1}_{j-1} = \mu^{k+1}_j N_j+ P^{k+1}_j A^{k+1}_j, \qquad  A^{k+1}_{j} := N_j - \Delta s \; \sum_{0 <s_i \leq \tau_j} \rho^{k+1}_{j,i} .
\end{equation}
As observed previously, these are $2L-1$ parameters for $L$ constraints and several options are possible.

\begin{table}[p] \centering
	\begin{tabular}{cccccc}
		\hline\hline
		$j  $  & 1   &2     & 3    & 4  &5  \\
		\hline
		$N_j$ & 8000 & 4000 & 2500 &1000 & 500 \\
		$\mu_j$& 1.6  &  1.6 &  1.6 &  1.6 & 5. \\
		$\tau_j$&4  &4 & 4 &4 &4  \\
		$P_j$& 3.37  &5. & 2.32 & 5. & --  \\
		$h_j$& 0.32  &0. & 0.05 & 0. & 0.1  \\
		$T_j$& 0.39 & .30 & .51& .30 & 2.1 \\
		$RP_j$& 0.05  & 0.03& 0.07 & 0.03 & 0.14
		\\
		\hline\hline
	\end{tabular}
	\\[-2mm]
	\caption{(High turnover) The company has to recruit externally at some job levels. The maximum internal hiring rate has been fixed to be $P_{\rm max}=5$. See Table 1  for notations.}
	\label{table2} \end{table}	
\paragraph{Maximize the internal promotion} To achieve this goal, one can choose to take the largest value $ P_{j-1}^{k+1}$ satisfying \eqref{ncchoice}, imposing however that $P_{j-1}^{k+1} \leq P_{\rm max} $ in order to avoid the subtle interpretation of $P_{j-1}^{k+1}=+\infty$. This yields the following formula:
$$
P_{j-1}^{k+1} = \min \left(P_{\rm max} , \f{ \mu^{k+1}_j N_j+ P^{k+1}_j A^{k+1}_j   }{ A^{k+1}_{j-1} }\right) .
$$
If $P_{j-1}^{k+1} < P_{\rm max}$, we choose $h_{j}^{k+1}=0$. If $P_{j-1}^{k+1} = P_{\rm max}$, we choose $h_{j-1}^{k+1}$ in order to achieve the equality~\eqref {ncchoice}.

\paragraph{Impose a proportion of external hiring} We may also choose to impose a given minimal proportion $\al$ of external hirin. This means $h_{j}^{k+1} = \al P_{j-1}^{k+1} + \delta_{j}^{k+1}$, with $\delta_{j-1}^{k+1} \geq 0$  nonnegative being used when the $\al$-strategy is not  enough to fill the job level $j-1$.  Then, we obtain the formula:
$$
P_{j-1}^{k+1} = \min \left(P_{\rm max} , \f{ \mu^{k+1}_j N_j+ P^{k+1}_j A^{k+1}_j   }{ (1+\al) A^{k+1}_{j-1} }\right) .
$$
If $P_{j-1}^{k+1} < P_{\rm max}$, we choose $h_{j}^{k+1}= \al P_{j-1}^{k+1}$. If $P_{j-1}^{k+1} = P_{\rm max}$, we choose $\delta_{j}^{k+1}$ in order to achieve the equality~\eqref {ncchoice}. 

\paragraph{Numerical examples} We have tested the algorithm with data describing an internal labor market with a “low” attrition rate. We concentrate ourselves on its steady state. The Table~\ref{table1} gives the data $N_j$, $\tau_j$ and the corresponding maximal internal promotion rate $P_j$. We also compute the proportion of employees ready to be promoted $RP_j= A_j/N_j$ and the average promotion when ready to be promoted
$$
T_j = \Delta s \sum_{s_i > \tau_j } s_i \rho^{k}_{j,i} .
$$
With low attrition, external hiring can be avoided.
\\

Table~\ref{table2} gives the results with a higher attrition rate and external hiring is needed. 

\subsection{Numerical simulation of a cost optimal structure with flexible workers}
\label{sec:numFlex}

Under the proposed organizational model described in section~\ref{sec:temporary}, it is possible to minimize a company labor cost footprint. The main question behind the minimization problem is about understanding the value of paying a flat premium to outsource some activities versus having to support constant cost increases due to a pay for tenure type of wage scheme.  As previously mentioned, this problem is difficult to solve and does not lead to closed formulas. In order to numerically approximate the cost optimal state of the organization, genetic algorithms (referred to as GAs) have been used. \\
GAs consist in search heuristics inspired by the basic principles of biological evolution \cite{GApackage}. In a nutshell, they start with a set of candidate solutions (population), and for each iteration (generation), the better solutions are selected (parents) and used to generate new solutions (offspring). The generation of new solutions can follow a variety of pattern: parent information can be recombined (crossover) or randomly modified (mutation) etc… The offspring are then used to update the population (evolution) and replace its weakest members. After each iteration, the overall population gets closer to a minimum. GAs have empirically proven useful to reach global minimum in a number of iterations that is quadratic in the number of sought parameters (see \cite{kumar2010genetic} for a review).
\\
In the case at hand, a generic GA implementation in  the software $R$ leveraging the $genalg$ package \cite{lucasius:1993:uugapcpc} was used. Once natural boundaries were selected for the parameters of interest (namely $\alpha_j \& p_j$), a population size of 200 individuals and a total of 250 iterations were fixed with the default implementation of the algorithm. Note that we specified a mutation chance of 10\%. The organization parameters that were used are described in Table~\ref{table3}  and echo the previous subsection's simulations. As per wage scheme, a 40\% increase per hierarchical level was assumed based on a 35\$/h salary at the lowest organizational level ($j=1$). For the sake of simplicity, the wage increase rate associated to tenure $r$ will be assumed at 4\%.
\begin{table}[h] \centering
	\begin{tabular}{cccccc}
		\hline\hline
		Org. Level $j$   & 1   &2     & 3    & 4  &5 \\
		\hline
		$N_j$   &5500 & 5200     &3800  &1800  &500 \\
		$\mu_j$ [\%] &8 & 8 & 8 & 8 & 20 \\	
		$\tau_j$ [Years]&4    &4        &4    &4     & 4 \\
		$w^0_j$ [\$/h] &35       &49    &69    &96   & 134 \\
		\hline\hline
	\end{tabular}
	\\[-2mm]
	\caption{Organization model parameters.}
	\label{table3} \end{table}
In the developed example, a uniform premium $B$ was assumed across all levels so that
$$
 w_j^t=(1+B).w_j^p, \qquad j=1,..., 5.
$$
We have tested three cases, $B=5\%$, $B=10\%$ and $B=20\% $ which we have  compared to the no-contractor scenario, i.e.,  $B=\infty$. Results are displayed in the  Table~\ref{table4}.
\begin{table}[h] \centering
	\begin{tabular}{c|cccccc|c}
		\hline\hline
		Contractor premium &Org. Level $j$   & 1   &2     & 3    & 4  & 5 &  Optimal Org. Cost [M\$/h]\\
		\hline
		$B=\infty$ & $\alpha_j$   & - &1.00     &1.00  &1.00  &1.00& 1,13\\
		\hline
		$B=20\%$& $\alpha_j$   & - &  1.32 & 1.04 & 1.02 & 1.00 & \\
		&$p_j$&0.23 &0.22 &0.22 &0.39& 0.99& 1,10\\
		\hline
		$B=10\%$& $\alpha_j$   & - &3.89& 1.02& 1.851& 1.29& \\
		&$p_j$&0.02 &0.075 &0.07 &0.20& 0.99 &1,04 \\
		\hline
		$B=5\%$& $\alpha_j$   & - &2.31 &1.08& 1.75& 1.42 &\\
		&$p_j$& 0.03& 0.06 &0.07 &0.18 &0.99& 0,99\\
		\hline\hline
	\end{tabular}
	\\[-2mm]
	\caption{Simulation results with temporary hirings.}
	\label{table4} \end{table}

The organizational footprint that is being discussed here can build upon itself internally  for a total of 1,13M\$/h, but  does not reach optimality if no flexibility is considered. Leveraging an internal labor market however leads to a certain degree of stagnation. As such leveraging flexible workforce solutions becomes financially interesting under the previous assumption for a premium $B<20\%$. The associated financial gain is of course increasing as the contractor premium $B$ gets lower. As such a 5\% premium yields more than 12\% in savings. Finally note that, in all our experiments, the top of the organization ($j=5$) never gets outsourced. However the proportion of external vs internal sourcing (represented by $\alpha$) increases at level 4 when flexible workforce solutions are cheap.


\section{Conclusion and perspectives}\label{sec:conc}

In this paper, a scalable hierarchical organization model was built. Its main features were to account for employee promotional readiness while controlling labor costs. The model shows that organizational inertia may come at a cost depending in its set up. Because of the organization needs and design, employees are indeed potentially waiting for a new internal opportunity while performing the same task. This occurs in situations with low turnover rates and  is detrimental for both the company that pay them to stay (henceforth increasing its operational costs) and for the individual.

The proposed solution investigates the notion of 
temporary labor. Rather than using flexibility to adapt to a volatile workload, the ideas developed in this paper are to optimally contract a part of this organization to create a lean and cost optimal environment for the individuals that belong to the company. The main results shows that it is possible to do so with various flexibility solutions (floaters, temporary workers) and with different model of employee commitment (employee grown internally vs brought from the outside).

If the proposed method is helpful in setting up an organization target, its operationalization would have to be discussed. The primary question that arises lies in the potential use of huge chunks of flexible workers.  Activities may not necessarily be available for outsourcing (ex: proprietary methods or processes).
Behind this potential core versus non-core activities discussion, the introduction of flexible workers may not call for the same type of workers in the “core” activities. For instance, for permanent people manager will have to deal with important workforce changes. While coaching their team members, they will have to adapt their efforts between permanent and flexible workers.
This organizational method also calls for a robust succession process especially in the promotion attributions that is not described in this paper. We believe this discussion should be skill based and requires the definition of company-wide standards. To this extent, note that a lot of ground has already been covered by principal-agent theoretical approaches (see \cite{ 10.2307/2564725} for example). 
Finally another potential operational issue is about performance and productivity. Creating change in the workforce to optimize costs may come at a performance loss if not managed correctly. This calls for caution in potential implementation as the resulting balance may not beneficial to the organization.


\bibliographystyle{ieeetr}  
\bibliography{OrgSpeedV11}

\begin{thebibliography}{10}

\bibitem{Bordoloi2001169}
S.~K. Bordoloi and H.~Matsuo, ``Human resource planning in knowledge-intensive
  operations: A model for learning with stochastic turnover,'' {\em European
  Journal of Operational Research}, vol.~130, no.~1, pp.~169 -- 189, 2001.

\bibitem{Burger_2016}
M.~Burger, A.~Lorz, and M.-T. Wolfram, ``On a {B}oltzmann mean field model for
  knowledge growth,'' {\em SIAM J. Appl. Math.}, vol.~76, no.~5,
  pp.~1799--1818, 2016.

\bibitem{VandenBergh2013367}
J.~V. den Bergh, J.~Beliën, P.~D. Bruecker, E.~Demeulemeester, and L.~D.
  Boeck, ``Personnel scheduling: A literature review,'' {\em European Journal
  of Operational Research}, vol.~226, no.~3, pp.~367 -- 385, 2013.

\bibitem{Klehe2011217}
U.-C. Klehe, J.~Zikic, A.~E.~V. Vianen, and I.~E.~D. Pater, ``Career
  adaptability, turnover and loyalty during organizational downsizing,'' {\em
  Journal of Vocational Behavior}, vol.~79, no.~1, pp.~217 -- 229, 2011.

\bibitem{balancing}
Komarudin, T.~De~Feyter, M.-A. Guerry, and G.~Vanden~Berghe, ``Balancing
  desirability and promotion steadiness in partially stochastic manpower
  planning systems,'' {\em Comm. Statist. Theory Methods}, vol.~45, no.~6,
  pp.~1805--1818, 2016.

\bibitem{henderson}
J.~S. Henderson, ``Stochastic optimal control of internal hierarchical labor
  markets,'' {\em J. Optim. Theory Appl.}, vol.~30, no.~1, pp.~99--115, 1980.

\bibitem{ribes:hal-01556746}
E.~Ribes, K.~Touahri, and B.~Perthame, ``{Employee turnover prediction and
  retention policies design: a case study}.'' working paper or preprint, 2017.

\bibitem{Qin201519}
R.~Qin, D.~A. Nembhard, and W.~L.~B. II, ``Workforce flexibility in operations
  management,'' {\em Surveys in Operations Research and Management Science},
  vol.~20, no.~1, pp.~19 -- 33, 2015.

\bibitem{Wild199395}
B.~Wild and C.~Schneewei, ``Manpower capacity planning - a hierarchical
  approach,'' {\em International Journal of Production Economics}, vol.~30-31,
  no.~C, pp.~95--106, 1993.

\bibitem{Harper2010}
P.~R. Harper, N.~H. Powell, and J.~E. Williams, ``Modelling the size and
  skill-mix of hospital nursing teams,'' {\em Journal of the Operational
  Research Society}, vol.~61, no.~5, pp.~768--779, 2010.

\bibitem{Perthame_book}
B.~Perthame, {\em Transport equations in biology}.
\newblock Frontiers in Mathematics, Birkh\"auser Verlag, Basel, 2007.

\bibitem{gurtin1974non}
M.~E. Gurtin and R.~C. MacCamy, ``Non-linear age-dependent population
  dynamics,'' {\em Archive for Rational Mechanics and Analysis}, vol.~54,
  no.~3, pp.~281--300, 1974.

\bibitem{farkas2002stability}
M.~Farkas, ``On the stability of stationary age distributions,'' {\em Applied
  mathematics and computation}, vol.~131, no.~1, pp.~107--123, 2002.

\bibitem{Doumic2017}
M.~Doumic, B.~Perthame, E.~Ribes, D.~Salort, and N.~Toubiana, ``Toward an
  integrated workforce planning framework using structured equations,'' {\em
  European Journal of Operational Research}, vol.~262, pp.~217--230, 2017.

\bibitem{bouchut2004nonlinear}
F.~Bouchut, {\em Nonlinear stability of finite volume methods for hyperbolic
  conservation laws and well-balanced schemes for sources}.
\newblock Series Frontiers in Mathematics, Birkh{\"a}user, Basel, 2004.

\bibitem{leveque1992numerical}
R.~J. LeVeque, {\em Numerical methods for conservation laws}.
\newblock Birkh{\"a}user, Verlag, Boston, 1992.

\bibitem{GApackage}
L.~Scrucca, ``{GA}: A package for genetic algorithms in {R},'' {\em Journal of
  Statistical Software}, vol.~53, no.~4, pp.~1--37, 2013.

\bibitem{kumar2010genetic}
M.~Kumar, M.~Husian, N.~Upreti, and D.~Gupta, ``Genetic algorithm: Review and
  application,'' {\em International Journal of Information Technology and
  Knowledge Management}, vol.~2, no.~2, pp.~451--454, 2010.

\bibitem{lucasius:1993:uugapcpc}
C.~B. {Lucasius, Jr.} and G.~Kateman, ``Understanding and using genetic
  algorithms. {P}art 1. concepts, properties and context,'' {\em Chemometrics
  and {I}ntelligent {L}aboratory Systems}, vol.~19, pp.~1--33, 1993.
\newblock reprint in: Carlos B. {Lucasius, Jr.}, Towards Genetic Algorithm
  Methodology in Chemometrics.

\bibitem{10.2307/2564725}
C.~Prendergast, ``The provision of incentives in firms,'' {\em Journal of
  Economic Literature}, vol.~37, no.~1, pp.~7--63, 1999.

\end{thebibliography}
\end{document}